\documentclass[11pt,a4paper]{amsart}
\usepackage[margin=25mm]{geometry}
\usepackage[numbers]{natbib}
\usepackage{hyperref}

\usepackage{amssymb,amsmath}
\usepackage{amsthm}
\usepackage{bbm, bm}
\usepackage{stmaryrd}
\usepackage{color}
\usepackage{graphicx}
\usepackage{caption}
\usepackage{float}
\usepackage{wrapfig}
\usepackage{orcidlink}
\usepackage{amsaddr}
\usepackage{cleveref}
\usepackage{enumerate}
\usepackage{lineno}
\usepackage{mathrsfs}

\usepackage{abstract}

\input xy
\xyoption{all} 

\numberwithin{equation}{section}

\newtheorem{theorem}{Theorem}[section]
\newtheorem{corollary}{Corollary}[section]
\newtheorem{lemma}{Lemma}[section]
\newtheorem{proposition}{Proposition}[section]

\theoremstyle{definition}
\newtheorem{definition}{Definition}[section]
\newtheorem{example}{Example}[section]

\newenvironment{warning}[1][Warning.]{\begin{trivlist}
\item[\hskip \labelsep {\bfseries #1}]}{\end{trivlist}}

\newenvironment{remark}[1][Remark.]{\begin{trivlist}
\item[\hskip \labelsep {\bfseries #1}]  }{ \end{trivlist}}

\newcommand{\Id}{\mathbbmss{1}}

\newcommand{\rmi}{ \textnormal{i}}
\newcommand{\rmd}{\textnormal{d}}

\DeclareMathOperator{\Vect}{Vect}

\DeclareMathOperator{\Span}{Span}

\newcommand{\catname}[1]{\textnormal{\texttt{#1}}}

\font\black=cmbx10 \font\sblack=cmbx7 \font\ssblack=cmbx5 \font\blackital=cmmib10  \skewchar\blackital='177
\font\sblackital=cmmib7 \skewchar\sblackital='177 \font\ssblackital=cmmib5 \skewchar\ssblackital='177
\font\sanss=cmss10 \font\ssanss=cmss8 
\font\sssanss=cmss8 scaled 600 \font\blackboard=msbm10 \font\sblackboard=msbm7 \font\ssblackboard=msbm5
\font\caligr=eusm10 \font\scaligr=eusm7 \font\sscaligr=eusm5  \font\fraktur=eufm10
\font\sfraktur=eufm7 \font\ssfraktur=eufm5 
\font\bsymb=cmsy10 scaled\magstep2
\def\all#1{\setbox0=\hbox{\lower1.5pt\hbox{\bsymb
       \char"38}}\setbox1=\hbox{$_{#1}$} \box0\lower2pt\box1\;}
\def\exi#1{\setbox0=\hbox{\lower1.5pt\hbox{\bsymb \char"39}}
       \setbox1=\hbox{$_{#1}$} \box0\lower2pt\box1\;}

\def\tx#1{{\fam0\relax#1}}

\newfam\bifam
\textfont\bifam=\blackital \scriptfont\bifam=\sblackital \scriptscriptfont\bifam=\ssblackital

\newfam\blfam
\textfont\blfam=\black \scriptfont\blfam=\sblack \scriptscriptfont\blfam=\ssblack

\newfam\bbfam
\textfont\bbfam=\blackboard \scriptfont\bbfam=\sblackboard \scriptscriptfont\bbfam=\ssblackboard

\newfam\ssfam
\textfont\ssfam=\sanss \scriptfont\ssfam=\ssanss \scriptscriptfont\ssfam=\sssanss
\def\sss#1{{\fam\ssfam\relax#1}}

\newfam\clfam
\textfont\clfam=\caligr \scriptfont\clfam=\scaligr \scriptscriptfont\clfam=\sscaligr

\newfam\frfam
\textfont\frfam=\fraktur \scriptfont\frfam=\sfraktur \scriptscriptfont\frfam=\ssfraktur

\def\hpb#1{\setbox0=\hbox{${#1}$}
    \copy0 \kern-\wd0 \kern.2pt \box0}
\def\vpb#1{\setbox0=\hbox{${#1}$}
    \copy0 \kern-\wd0 \raise.08pt \box0}

\def\pmb#1{\setbox0\hbox{${#1}$} \copy0 \kern-\wd0 \kern.2pt \box0}
\def\pmbb#1{\setbox0\hbox{${#1}$} \copy0 \kern-\wd0
      \kern.2pt \copy0 \kern-\wd0 \kern.2pt \box0}
\def\pmbbb#1{\setbox0\hbox{${#1}$} \copy0 \kern-\wd0
      \kern.2pt \copy0 \kern-\wd0 \kern.2pt
    \copy0 \kern-\wd0 \kern.2pt \box0}
\def\pmxb#1{\setbox0\hbox{${#1}$} \copy0 \kern-\wd0
      \kern.2pt \copy0 \kern-\wd0 \kern.2pt
      \copy0 \kern-\wd0 \kern.2pt \copy0 \kern-\wd0 \kern.2pt \box0}
\def\pmxbb#1{\setbox0\hbox{${#1}$} \copy0 \kern-\wd0 \kern.2pt
      \copy0 \kern-\wd0 \kern.2pt
      \copy0 \kern-\wd0 \kern.2pt \copy0 \kern-\wd0 \kern.2pt
      \copy0 \kern-\wd0 \kern.2pt \box0}

\mathchardef\za="710B  
\mathchardef\zb="710C  
\mathchardef\zg="710D  
\mathchardef\zd="710E  
\mathchardef\zve="710F 
\mathchardef\zz="7110  
\mathchardef\zh="7111  
\mathchardef\zvy="7112 
\mathchardef\zi="7113  
\mathchardef\zk="7114  
\mathchardef\zl="7115  
\mathchardef\zm="7116  
\mathchardef\zn="7117  
\mathchardef\zx="7118  
\mathchardef\zp="7119  
\mathchardef\zr="711A  
\mathchardef\zs="711B  
\mathchardef\zt="711C  
\mathchardef\zu="711D  
\mathchardef\zvf="711E 
\mathchardef\zq="711F  
\mathchardef\zc="7120  
\mathchardef\zw="7121  
\mathchardef\ze="7122  
\mathchardef\zy="7123  
\mathchardef\zf="7124  
\mathchardef\zvr="7125 
\mathchardef\zvs="7126 
\mathchardef\zf="7127  
\mathchardef\zG="7000  
\mathchardef\zD="7001  
\mathchardef\zY="7002  
\mathchardef\zL="7003  
\mathchardef\zX="7004  
\mathchardef\zP="7005  
\mathchardef\zS="7006  
\mathchardef\zU="7007  
\mathchardef\zF="7008  
\mathchardef\zW="700A  
\mathchardef\zC="7009  

\newcommand{\be}{\begin{equation}}
\newcommand{\ee}{\end{equation}}

\newcommand{\bea}{\begin{eqnarray}}
\newcommand{\eea}{\end{eqnarray}}
\def\*{{\textstyle *}}
\newcommand{\R}{{\mathbb R}}

\newcommand{\C}{{\mathbb C}}
\newcommand{\Z}{{\mathbb Z}}

\newcommand{\s}{{\textstyle *}}

\def\ul{\underline}
\def\Sec{\sss{Sec}}
\def\Vect{\sss{Vect}}

\def\sT{{\sss T}}

\def\xi{\tx{i}}

\def\s*{{\scriptstyle *}}

\def\ul{\underline}

\newcommand{\beas}{\begin{eqnarray*}}
\newcommand{\eeas}{\end{eqnarray*}}

\title{Biunital Lie Semiheaps and  Lie  Triple Systems} 
\author{Andrew James Bruce \,\orcidlink{0000-0001-8197-2263}} 
  \address{Independent Researcher\\ Cardiff, United Kingdom  }  
   \email{andrewjamesbruce@googlemail.com}

   \date{\today}
 
\begin{document}
 \maketitle
\vspace{-20pt}
\begin{abstract}{\noindent  We examine the differential geometry of Lie semiheaps that admit biunital elements. Open Ol'shanskii semigroups, which appear widely across causal geometry and harmonic analysis, provide a large class of examples. It is shown that the subset of biunits is a smoothly embedded Lie subheap. Moreover, there is a canonical Lie triple algebroid associated with any proper biunital  Lie semiheap. That is, there is a vector bundle over the manifold of biunits that comes equipped with a ternary bracket that satisfies the axioms of a  Lie triple system together with a compatible `higher' anchor map. Although every Lie triple algebroid associated with a proper biunital Lie semiheap is trivialisable, it is not canonically so, unlike the case of the tangent bundle of a Lie group.  As concrete examples, we examine the Lie semiheap of real $n\times n$ matrices, and the conformal causal Lie semiheap associated with the Lie group $\mathrm{SU}(2,2)$, which is a double cover of the conformal group of Minkowski spacetime, and the future-directed light cone. }\\

\noindent {\Small \textbf{Keywords:}~Lie Semiheaps;~Biunits;~Lie Triple Systems;~Lie Triple  Algebroids;~Non-Invertible Transformations}
\end{abstract}
\tableofcontents
\section{Introduction and Motivation}
\subsection{Introduction and Background} 
Lie groups, and more generally, Lie groupoids play a fundamental role in differential geometry and geometric approaches to physics by providing a robust notion of continuous transformations and symmetries.  Heaps and semiheaps have a long algebraic history going back to the 1920's with Baer  \cite{Baer:1929} and Pr\"ufer \cite{Prufer:1924}. Heuristically, we think of a heap as a group for which the identity element has been forgotten; a little more carefully, we have a ternary operation together with a generalised version of associativity.  Indeed, given a group $(G, \cdot)$, we can construct a heap by defining $[g_1, g_2, g_3] := g_1 \cdot g^{-1}_2 \cdot g_3$.  In fact, up to isomorphism, every heap arises from a group in this way.   Conversely, given a heap $(H, [-,-,-])$, by selecting any point $h \in H$, we can construct a group by defining $h_1 \cdot_h h_2 :=  [h_1, h, h_2]$. Note that the element $h$ is the identity element and that inversion is given by  $h_1^{-1} = [h, h_1 ,h]$. However, we do not have an equivalence of categories here as passing from a heap to a group is not natural; morphisms of heaps need not respect the selected points. Semiheaps, as introduced by Wagner\footnote{In translations of Soviet era journals, Wagner is spelt Vagner.} \cite{Vagner:1953}, are a slightly more general notion where the association with groups is lacking. However, it is well-known that any semiheap can be embedded into an involutive monoid.\par
Lie heaps and semiheaps were introduced by the author in \cite{Bruce:2024}\footnote{Maybe better stated as re-introduced, as Konstantinova \cite{Konstantinova:1978} discusses semiheap bundles. However, this work is in Russian, and the author has not been able to obtain a copy or translation.} as heap and semiheap objects in the category of smooth manifolds. In plain language, we have analogues of Lie groups for heaps and semiheaps. The main complication working with general Lie semiheaps when compared with the standard situation of Lie groups is that we do not have an identity element; in effect, there may be many identities or none at all,  and left/right translations are not diffeomorphisms. Conceptually, the absence of a globally fixed identity should remind the reader of a principal $G$-bundle, which can, heuristically, be viewed as a smoothly parameterised family of Lie groups (the fibres) for which the unit in each fibre is not fixed. The selection of a unit in each fibre requires a local trivialisation to be defined. However, whereas principal bundles rely on group actions and strictly invertible fibre transformations, Lie semiheaps naturally handle non-invertible transformations.\par 
We will concentrate on Lie semiheaps $(S, [-,-,-])$ that are not Lie heaps, but contain biunits, that is, elements $e \in S$, such that $[e,e,s] = [s,e,e] = s$ for all $s \in S$. If all the elements are biunits, then we have a Lie heap.  Such Lie semiheaps we will refer to as \emph{proper biunital Lie semiheaps}.   By working with this particular class of Lie semiheaps, we can locally apply notions and theorems from the theory of Lie groups and algebras without the geometry globally collapsing to that of Lie groups and algebras. Philosophically, we think of proper biunital semiheaps as Lie monoids for which the identity has been `forgotten' or not  `fixed'. We remark that monoid actions on smooth manifolds are not well-studied; the interested reader should consult \cite{Jozwikowski:2016} and references therein. In particular, the multiplicative monoid of reals is fundamental in (non-negatively) graded geometry, see Grabowski \& Rotkiewicz \cite{Grabowski:2012}. \par 
\smallskip 

\noindent\textbf{Key Results:} We establish two fundamental theorems about the differential geometry of proper biunital Lie semiheaps. 
\begin{description}
\item[Theorem \ref{thm:BiunitEmbed}] The set of biunits is a smoothly embedded Lie subheap of the proper biunital Lie semiheap.
\item[Theorem \ref{the:LieTripSys}] Associated with any proper biunital Lie semiheap is a vector bundle over the manifold of biunits whose sections carry the structure of a Lie triple algebroid. 
\end{description}
Lie triple systems owe their inception to Cartan \cite{Cartan:1926}  and his classification of symmetric spaces. The algebraic theory is due to  Jacobson \cite{Jacobson:1949}, and Yamaguti \cite{Yamaguti:1960} constructed the representation and cohomology theory of Lie triple systems. Remarkably, symmetric spaces and Lie triple systems have found application in the numerical analysis of differential equations, see Munthe-Kaas et al. \cite{Munthe-Kaas:2014}.  We stress that a Lie triple system is defined on a single vector space, whereas Lie triple algebroids are defined on the module of sections of a vector bundle (see Definition \ref{def:LieTripSymAlg}). Lie triple systems associated with symmetric spaces rely on the existence of an involution acting on the tangent spaces or the associated Lie algebra of the symmetric space and its associated eigenspace decompositions. The same core mechanism is behind the Lie triple algebroids associated with proper biunital Lie semiheap, though the origin of the involution is, in general, different. In particular, the differential of the involution $\mathsf{I}_e(-):= [e,-,e]$ for every biunit $e$ provides the involution on the tangent space at each biunit. We highlight that the underlying vector bundle structure of a Lie triple algebroid associated with a proper biunital Lie semiheap is trivialisable; however, not canonically.  The selection of a biunit is required to `fix' the trivialisation. Heuristically, we interpret the selection of a biunit and the resulting trivialisation as picking an observer and using their frame to describe the Lie triple algebroid.   \par 
As a particular class of proper biunital Lie semiheaps, we have open Ol'shanskii semigroups (see Example \ref{ema:OlSemi}), as found in modern Lie theory, causal geometry, constructive field theory, and harmonic analysis (see \cite{Hilgert:1989}). We explore in some detail the proper biunital Lie semiheap of $n\times n$ real matrices (see Example \ref{exa:nxnMatrices}).  A physically relevant example is the conformal causal Lie semiheap associated with the Lie group $\mathrm{SU}(2,2)$ and the future-directed light cone (see Example \ref{exa:ConfCausal}). In the context of causal conformal mechanics, biunits here are identified with complete observers, thus realising the desire to formulate aspects of physics without a choice of observers. 
%
%
\subsection{Motivation}
There has been a renewed interest in ternary structures in algebra and geometry. Generically, one can think of geometric structures based on heaps and closely related ternary operations as `affine geometric structures' in the sense that some privileged element is not fixed. Thus, developing ternary structures is seen as part of developing structures that do not rely on an underlying vector space but are, in some sense, affine.  Selecting a privileged element to `anchor' the structure is analogous to choosing an observer or picking a gauge in physics.  If one wants to develop physics without reference to observers, heaps and related ternary structures potentially offer a mathematical route.  As an example, Ostermann \& Schmidt \cite{Ostermann:1966} showed that affine spaces can be defined without reference to vector spaces by using two ternary operations; one entirely on the set and the other representing an action of the field of scalars. This was expounded upon by Breaz et al. \cite{Breaz:2024}. Brzeziński \cite{Brzezinski:2022}  reformulated Lie brackets on affine spaces/bundles using ternary structures and so produced a frame-independent understanding of Lie affebras as introduced by Grabowska et al. \cite{Grabowska:2003}. The  affinisation of linear algebra structures, and in particular Leibniz algebras, was presented by Brzeziński et al. \cite{Brzezinski:2026}. \par 
Focussing in on proper biunital Lie semiheaps, it is important to note that the left and right translations (see Definition \ref{def:RLActions}) form a monoid rather than a group; we have smooth maps that are not invertible. One might ask, ``Why not start with a Lie monoid?''.  The answer is that a Lie monoid has a single identity which we may heuristically interpret as a privileged  observer. As we desire a geometric formalism that allows for many ``identities'' or ``observers'', we consider Lie semiheaps.  The infinitesimal theory of proper biunital Lie semiheaps is geometrically richer than that of Lie monoids, as we do not, in general, have a single biunit.  Thus, Lie semiheaps offer a framework for smooth and infinitesimal non-invertible symmetries in classical physics and geometric mechanics beyond Lie groups/monoids.
%
%
\subsection{Comparison with Lie Groupoids}
To conceptually aid the reader, we provide a short table comparing proper biunital Lie semiheaps with Lie groupoids. For the (abstract) definition of a  Lie triple algebroid, see Definition \ref{def:LieTripSymAlg}. Our standard reference for Lie groupoids and Lie algebroids is Mackenzie \cite{Mackenzie:2005}. Indeed, we consider proper biunital Lie semiheaps as `cousins' of Lie groupoids.  Similarly, Lie triple algebroids should be viewed as a `ternary relative' of Lie algebroids. 
\begin{table}[H]
\centering
{\renewcommand{\arraystretch}{1.2}
\begin{tabular}{|l|p{5cm}|p{5cm}|}
\hline
 & \textbf{Proper Biunital Lie Semiheaps} & \textbf{Lie Groupoids} \\ \hline
\textbf{Object} & A proper semiheap object in the category of smooth manifolds that possesses biunits. & A groupoid object in the category of smooth manifolds. \\ \hline
\textbf{Structure} & Ternary composition without an inverse map. &  Binary partial-composition with an inverse map. \\ \hline
\textbf{Identities} & A smoothly embedded Lie subheap of biunits. & A smoothly embedded submanifold of units. \\ \hline
\textbf{Algebroid} & \textbf{Lie triple algebroid}: A vector bundle over the manifold of biunits whose sections form a Lie triple system. & \textbf{Lie algebroid}: A vector bundle over the manifold of units whose sections form a Lie algebra. \\ \hline
\textbf{Anchor Map} & Binary anchor that is skew-symmetric together with a generalised Leibniz rule. & Unary anchor together with the Leibniz rule. \\ \hline
\end{tabular}}
\end{table} 
\noindent \textbf{Comment.}  The related notions of Lie triple system bundles and Lie triple system algebroids have independently been proposed by Banerjee \&  Mukherjee \cite{Banerjee:2026}.  \par 
Note that the subset of units (identities/objects) of a Lie groupoid is by definition a smoothly embedded submanifold. For the case of proper biunital Lie semiheaps, the embedding of the biunits is automatic and not imposed by definition.\par 
We remind the reader that a principal $G$-bundle $P \rightarrow M$ can be considered as a Lie groupoid by passing to the \emph{Atiyah groupoid} $P \times_G P \rightarrow M$, where we then identify the submanifold of units (identities) with $M$.
\begin{remark}
Partial (Lie) semiheaps were employed in \cite{Bruce:2023}  to study nowhere-vanishing null vector fields in general relativity.
\end{remark}
\smallskip 

\noindent\textbf{Organisation:} In Section \ref{subsec:RecLieSemi} we recall the notion of a Lie semiheap and the basic concepts we will need throughout this work. In Section \ref{subsec:BiuLieSemi},  we define proper biunital Lie semiheaps and establish the embedding of the biunits into the Lie semiheap. In Section \ref{subsec:InfTheBiuLieSemiHeaps}, we explore the infinitesimal theory of proper biunital Lie semiheaps and introduce the associated Lie triple algebroid.   We end with concluding remarks and suggestions for future work in Section \ref{subsec:ConRem}. The abstract definition of a Lie triple algebroid is presented in Appendix \ref{App:LieTripSymAlg}.
%
%
\section{Biunital Lie Semiheaps  \& Triple Systems}
\subsection{Recollection of Lie Semiheaps}\label{subsec:RecLieSemi}
Recall that a Lie group is a group object in the category of smooth manifolds, i.e., a smooth manifold equipped with three smooth maps, the unit, inverse and multiplication maps, that satisfy the standard axioms of a group (see \cite[III.6]{MacLane:1988} for group objects in a category).   A Lie semiheap is a semiheap object in the category of smooth manifolds. That is, a smooth manifold together with a ternary operation that satisfies the axioms of a semiheap; our main reference for heaps and semiheaps is \cite{Hollings:2017}.  
\begin{definition}\label{def:SemiHeap}
A set (possibly empty) equipped with a ternary product $(S, \,[-,-,-])$ is a \emph{semiheap} if the ternary product is para-associative\footnote{Para-associativity is also referred to as \emph{B-associativity} following Kerner \cite{Kerner:2018}. \emph{A-associativity} is defined as  $[s_1, s_2 , [s_3, s_4, s_5]] = [s_1, [s_2, s_3, s_4], s_5] = [[s_1, s_2, s_3], s_4, s_5]$. }, i.e., 
$$[s_1, s_2 , [s_3, s_4, s_5]] = [s_1, [s_4, s_3, s_2], s_5] = [[s_1, s_2, s_3], s_4, s_5]\, ,$$
for all $s_i \in S$. A semiheap is said to be an \emph{abelian semiheap} if $[s_1, s_2, s_3] = [s_3, s_2, s_1]$ for all $s_i \in S$.  A semiheap is a \emph{heap} if every element of the semiheap is \emph{biunitary}, i.e., the Mal'cev identity $s = [s',s',s] = [s,s',s' ]$ holds for all $s,s' \in S$.
A semiheap is said to be \emph{totally commutative} if all permutations of the entries are equal, i.e., $[s_1, s_2, s_3] = [s_2, s_1, s_3] = [s_1, s_3, s_2] = [s_3, s_2, s_1] = [s_2, s_3, s_1] = [s_3, s_1, s_2]$.
\end{definition}
\begin{definition}
Let $(S, \, [-,-,-])$ and $(S', \, [-,-,-]')$ be semiheaps. A \emph{homomorphism of semiheaps} is a map $\phi : S \rightarrow S' $ that satisfies
$$\phi[s_1, s_2, s_3] =  [\phi(s_1),  \phi(s_2), \phi(s_3)]'\,,$$
for all $s_i \in S$.
\end{definition} 
 We thus obtain the category of semiheaps, which we denote as $\catname{SHeap}$.
\begin{example}
Let $(\mathcal{S}, \cdot )$  be a commutative semigroup. Then $[x,y,z] := x\cdot y \cdot z$ defines a semiheap structure on the set $\mathcal{S}$. 
\end{example}
\begin{example}
 Let $(M,g)$  be a (pseudo-)Riemannian manifold, then $\Vect(M)$ can be equipped with the semiheap structure $[X,Y, Z] :=  g(X,Y)\, Z$. One can show that the null vector fields are closed under this ternary operation and so form a subsemiheap of $\big(\Vect(M), [-,-,-]\big)$. 
\end{example}
\begin{example}
Let $\mathcal{C}(E)$ be the affine space of linear connections on a vector bundle $\pi : E \rightarrow M$.  The ternary operation $[\nabla_1, \nabla_2, \nabla_3]:=  \nabla_1 - \nabla_2 + \nabla_3$  defines a heap structure on $\mathcal{C}(E)$. See \cite{Bruce:2024b} for more details.
\end{example}
\begin{example}
Given a $*$-algebra $(\mathcal{A}, \cdot , *)$, then $[a,b,c ] := a \cdot b^* \cdot c$ defines a semiheap structure. 
\end{example}
\begin{definition}[\cite{Bruce:2024}]
A \emph{Lie semiheap} is a semiheap object in the category of real, finite-dimensional, Hausdorff and second-countable smooth manifolds.  A \emph{Lie semiheap homomorphism} is a smooth map $\psi : S \rightarrow S'$ that is also a semiheap homomorphism.
\end{definition}
The category of Lie semiheaps we denote as $\catname{LieSHp}$. We have the full subcategory of Lie heaps, which we denote as $\catname{LieHp}$. The obvious forgetful functor $\textrm{F}: \catname{LieHp} \rightarrow \catname{LieSHp}$, which forgets the biunitary property of all elements, is a fully faithful embedding.  Thus, we view Lie heaps as Lie semiheaps with further properties. Furthermore,  there are other full subcategories  $\catname{AbLieSHp}$ of abelian Lie semiheaps, and $\catname{AbLieHp}$ of abelian Lie heaps. We will, when there is no chance of confusion, denote a Lie semiheap by just its underlying manifold $S$.
\begin{definition}
A Lie semiheap $S$ is said to be a \emph{proper Lie semiheap} if (as an object) it is not in the essential image of the forgetful functor $\textrm{F}: \catname{LieHp} \rightarrow \catname{LieSHp}$.
\end{definition}
In simpler language, a proper Lie semiheap is a Lie semiheap that is not a heap, i.e., not all the elements satisfy the Mal'cev identity. 
\begin{definition}
Let $S$  be a Lie semiheap. A \emph{Lie subsemiheap} is a subset $S' \subset S$ that is closed under the semiheap operation and is a smoothly embedded submanifold $\iota : S' \hookrightarrow S$.  If such an $S'$ is a heap, then we speak of a \emph{Lie subheap}.
\end{definition}
\begin{definition}\label{def:RLActions}
Let $S$ be a Lie semiheap and fix a pair $(s_1, s_2) \in S\times S$. The smooth map 
\begin{align*}
R_{s_1 s_2}: &~ S \longrightarrow S \,, && s \mapsto [s, s_1, s_2]
\end{align*}
is a \emph{right translation}. Similarly, the smooth map
\begin{align*}
L_{s_1 s_2}: &~ S \longrightarrow S\,,  && s \mapsto [ s_1, s_2, s]
\end{align*}
is a \emph{left translation}.
\end{definition}
We will denote the set of right translations of a Lie semiheap by $R(S)$ and the set of left translations of a semiheap by $L(S)$.
\begin{proposition}\label{prop:TransSemi}
Let $S$ be a Lie semiheap. Then the sets $R(S)$ and $L(S)$ of right and left translations, respectively, are semigroups. 
\end{proposition}
\begin{proof}
Follows via direct calculation (see \cite{Bruce:2024} for details).
\end{proof}
In particular, 
\begin{align}\label{eqn:RLSemiGrp}
&R_{s_3 s_4}\circ R_{s_1 s_2} = R_{s_1[s_2, s_3, s_4] }\,, && L_{s_1 s_2} \circ L_{s_3 s_4} = L_{[s_1, s_2, s_3] s_4}\,,
\end{align}
for all $s_i \in S$.
\begin{example}\label{exa:R}
With respect to the standard smooth structure on $\R$, we have an abelian Lie heap structure given by $[s_1, s_2, s_3]:= s_1 -s_2 + s_3$.
\end{example}
\begin{example}\label{exa:R+}
With respect to the standard smooth structure on $\R_{>0}$, we have an abelian Lie semiheap structure  given by $[s_1, s_2, s_3]:= s_1 + s_2 + s_3$. Note that we have a proper Lie semiheap, i.e., we do not have a heap.
\end{example}
\begin{example}\label{exa:Grp}
Given a Lie group $G$, the ternary operation $[g_1, g_2, g_3] := g_1 g^{-1}_2 g_3$ defines a Lie heap structure on $G$. All Lie heaps are of this form.
\end{example}
\begin{remark}
Lie superheaps were defined using the functor of points by the author in \cite{Bruce:2025}.
\end{remark}
%
%
\subsection{Biunital Lie Semiheaps \& the Embedded Submanifold of Biunits}\label{subsec:BiuLieSemi}
We now turn our attention to Lie semiheaps that contain biunits. This extra structure allows for some aspects of standard Lie theory to be locally applied and allows for richer geometric structures. The algebraic tools employed are present in Wagner \cite{Vagner:1953}; it is their application and adaptation to the smooth setting that is novel.   We direct the reader to Hollings \& Lawson \cite{Hollings:2017} for translations of Wagner's fundamental work on semiheaps
and related structures.
\begin{definition}
Let $S$ be a Lie semiheap. The subset of \emph{left biunits} and  the subset of \emph{right biunits} are defined as 
$$B_L := \big \{e_L \in S ~~|~~ [e_L,e_L, s] = s~~\text{for all}~ s \in S \big\}\,, \qquad B_R := \big \{e_R \in S ~~|~~ [s,e_R, e_R] = s~~\text{for all}~ s \in S \big\}\,,$$
respectively.  The subset of \emph{biunits} is defined as $B :=  B_L \cap B_R$.
\end{definition}
Note that these sets may be empty, the entire Lie semiheap, or proper subsets.
\begin{remark}
Biunits are also referred to as \emph{Mal’cev elements} in the literature. 
\end{remark}
\begin{definition}
A Lie semiheap $S$ is said to be a \emph{left (right) biunital Lie semiheap} if $B_L \neq \emptyset$ ($B_R \neq \emptyset$). A Lie semiheap $S$ is said to be a \emph{biunital Lie semiheap} if $B := B_L \cap B_R \neq \emptyset$. A \emph{proper biunital Lie semiheap} is a biunital Lie semiheap that is also a proper Lie semiheap.
\end{definition}
Note that if $B = S$, then $S$ is a Lie heap. 
\begin{remark}
In \cite{Bruce:2024}, the author defined biunital Lie semiheaps as pointed Lie semiheaps for which the privileged element is a biunit. In the current text, the definition is broader and only requires the set of biunits $B$ to be non-empty.
\end{remark}
\begin{definition}
Let $S$ and $S'$  be proper biunital Lie semiheaps. \emph{A homomorphism of proper biunital Lie semiheaps}  $\psi : S \rightarrow S'$ is a homomorphism of Lie semiheaps that satisfies $\psi(B) \subseteq B'$.  
\end{definition}
The resulting category we denote as $\catname{PBLieSHp}$.
\begin{example}
With respect to its standard smooth structure, $\R$ equipped with the semiheap structure $[s_1, s_2, s_3] = s_1 + s_2 + s_3$ is a proper biunital Lie semiheap as $B_L = B_R =\{ 0\}$. 
\end{example}
\begin{example}
Consider the open cylinder $S = S^1 \times \R_{> 0}$ equipped with its standard smooth structure.  We define the semiheap structure as 
$$[(\theta_1, r_1), (\theta_2, r_2), (\theta_3, r_3) ] = (\theta_1 - \theta_2 + \theta_2, r_1 r_2 r_3) \,.$$
We can view this example as the Cartesian product of two Lie semiheaps; the Cartesian product was shown to be the categorical product in \cite[Proposition 2.34]{Bruce:2024}. The set of biunits is clearly $B = S^1 \times \{ 1\} \cong S^1$. Thus, we have a proper biunital Lie semiheap.
\end{example}
\begin{example}
Consider the quaternions $\mathbb{H} \cong \R^4$ equipped with its standard smooth structure. Then,  using the standard conjugation of quaternions, we can define the smooth semiheap operation $[q_1 , q_2 , q_3] := q_1 q_2^* q_3$.   The set of biunits is $B \cong S^2$ (unit quaternions), which as a Lie group is diffeomorphic to $\mathrm{Sp}(1) \cong \mathrm{SU}(2)$.  
\end{example}
\begin{remark}
The previous example generalises to the Cartesian product of a Lie heap (which, of course, corresponds to a Lie group) and a proper biunital Lie semiheap. Moreover, the Cartesian product of two (or more) proper biunital Lie semiheaps is again a proper biunital Lie semiheap.  
\end{remark}
\begin{example}
Let $G$ be a Lie group equipped with a smooth involution $\Theta : G \rightarrow G$.  Defining $g* :=\Theta(g_2)^{-1}$, we have the Lie semiheap structure on $G$ given by  $[g_1, g_2, g_3] := g_1 g_2^*g_3$. Observe that $B = G^\Theta := \{g \in G~~|~~ \Theta(g) =g\}$, and so we obtain a proper biunital Lie semiheap provided $B \neq \emptyset$ and $B \neq G$.
\end{example}
\begin{example} [Open Ol'shanskii semigroups \cite{Olshankskii}] \label{ema:OlSemi}Let $G$ be a connected semisimple Hermitian Lie group with Lie algebra $\mathfrak{g}$. Consider its Cartan decomposition $\mathfrak{g} = \mathfrak{k}\oplus \mathfrak{p}$, where $\mathfrak{k}$ is the Lie algebra of a maximal compact subgroup $K \subset G$, and $\mathfrak{p}$ is the non-compact complement. Let $C \subset \mathfrak{p}$ be a chosen open, convex, pointed cone that is invariant under the adjoint action of $K$, i.e., $\mathrm{Ad}(K) C = C$. Pointed here means that $C \cap (-C) = \emptyset$. We denote the complexification of $G$  as $G_\C$. The open Ol'shanskii semigroup $S^\textnormal{o}$ is defined as the image of $G \times C \rightarrow G_\C$ of the map $(g, v) \mapsto g \, \exp(\rmi\, v)$. The involution on $S^\textnormal{o}$ is the restriction of the canonical involutive antiautomorphism of the complex Lie group $G_\C$.  Using the polar decomposition, $s = g \, \exp(\rmi\, v)$, the involution is $s^* = \exp(- \rmi\, v)\, \Theta(g)^{-1}$, where $\Theta$ is Cartan's involution. Recall that $K= G^\Theta := \{g \in G~~|~~ \Theta(g) =g\}$. We then define the Lie semiheap operation $[s_1, s_2, s_3] := s_1 s_2^* s_3$, noting that the set of biunits  $B \stackrel{\sim}{\rightarrow} K \times C$ using the polar decomposition. In the classical language of involutive semigroups, the binits are unitary elements. Thus, we obtain a proper biunital Lie semiheap. See \cite{Hilgert:1989} for details of Ol'shanskii semigroups.
\end{example}
\begin{remark}
Neeb \cite{Neeb:1998} established that open Ol'shanskii semigroups are examples of Stein manifolds.
\end{remark}
\begin{lemma}\label{lem:BLClosed}
Let $S$  be a Lie semiheap and let $B_L$ and $B_R$ be its subset of left biunits and subset of right biunits, respectively. Then both $B_L$ and $B_R$ are closed subsets of $S$. 
\end{lemma}
\begin{proof}
We will establish the lemma for $B_L$ as the case of $B_R$ follows via near identical arguments.   As  $S$ is Hausdorff, the diagonal $\Delta_S \subset S \times S$ is closed.  For a fixed point $s \in S$, the evaluation map $h_s : S \rightarrow S \times S$ given by $s' \mapsto ([s',s', s], s)$ is smooth as it is built from smooth operations.  Moreover, the set 
$$E_s = \big \{ e \in S ~~|~~[e,e,s] =s \big\} = h^{-1}_s(\Delta_S)\,,$$ 
is closed in $S$. Thus, $B_L = \bigcap_{s \in S} E_s$ is an intersection of closed sets, and so is itself closed in $S$. 
\end{proof}
\begin{corollary}\label{cor:BClosedSubSet}
Let $S$  be a Lie semiheap; then the subset of biunits $B := B_L \cap B_R$ is a closed subset in $S$.
\end{corollary}
\begin{lemma}\label{lem:BClosed}
Let $S$ be a Lie semiheap; then  $B := B_L \cap B_R $ is closed under the ternary operation. That is, if $e_1, e_2 , e_3 \in B$ then $[e_1, e_2, e_3] \in B$.  Moreover, as every element is biunital, $B$ is a heap.
\end{lemma}
\begin{proof}
Let $e_1, e_2, e_3 \in B$ and $s \in S$. Then
\begin{align*}
[[e_1, e_2, e_3], [e_1, e_2, e_3], s] & = [[[e_1, e_2, e_3],e_3, e_2], e_1, s]=[[e_1,e_2, [e_3, e_3,e_2]], e_1, s]\\
&=[[e_1, e_2, e_2], e_1 , s ]= [e_1, e_1, s]=s\,.
\end{align*}
Thus, $[e_1, e_2, e_3] \in B_L$. The calculation establishing that $[e_1, e_2, e_3] \in B_R$ is almost identical to the above and so is omitted.  As every element in $B$ is a biunit and $B$ is closed under the ternary semiheap operation, we do indeed obtain a Lie heap.
\end{proof}
\begin{lemma}\label{lem:RLDiffs}
Let $S$ be a biunital Lie semiheap. Then the left and right translations generated by biunits form groups.
\end{lemma}
\begin{proof}
For any pair $e_1, e_2 \in B$ we have the associated right translation map $R_{e_1 e_2}(-):=[-, e_1, e_2]$. Similarly we have the associated left translation map $L_{e_1 e_2}(-):= [e_1,e_2, -]$.  The inverses are $R^{-1}_{e_1 e_2} = R_{e_2 e_1}$ and  $L^{-1}_{e_1 e_2} = L_{e_2 e_1}$.   To see this, observe that 
\begin{align*}
& R_{e_1 e_2 }\big( R_{e_2 e_1}(s)\big) = [[s, e_1, e_2], e_2, e_1] = [s, [e_2,e_2, e_1], e_1]= [s, e_1,e_1]=s\,,\\
& L_{e_1 e_2}\big( L_{e_2 e_1}(s)\big) = [e_1, e_2 , [e_2, e_1, s]] = [e_1, [e_1, e_2, e_2], s] =[e_1,e_1, s] = s\,.
\end{align*}
Interchanging $e_1 \leftrightarrow e_2$, we establish that we have left and right inverses (as expected). Thus, as the right and left translations are smooth maps that have smooth inverses, they are diffeomorphisms. In other words, the right/left translations generated by biunits form a group. 
\end{proof}
\begin{lemma}\label{lem:BTrans}
Let $S$ be a biunital Lie semiheap. Then the groups of left and right translations generated by biunits act transitively on the subset of biunits $B \subseteq S$.
\end{lemma}
\begin{proof}
As $R_{e e_1}(e)= [e, e, e_1] =e_1$ and $L_{e_1 e}(e) = [e_1, e, e]= e_1$ for arbitrary $e\in B$ and any chosen $e_1 \in B$, the group actions of right and right translations generated by biunits restricted to $B$ are transitive.        
\end{proof}
\begin{proposition}\label{prop:BInvIm}
Let $S$ be a biunital Lie semiheap,  i.e., $B \neq \emptyset$. For an arbitrary point $e \in B$ we define a  map $\eta_e : S \rightarrow S \times S$ given by $\eta_e(s) = \big( [s,s, e],[e,s,s]\big)$. Then $B = \eta^{-1}_e(e,e)$. Moreover, $\big(R_{e_1 e} \times L_{e e_1} \big)\big(\eta_{e_1}(s) \big) = \big([[s,s,e_1], e_1,e], [e , e_1,[e_1 , s,s ]] \big)  = \eta_e(s)$ for all $e_1 \in B$.
\end{proposition}
\begin{proof}
If $s \in B$, then $\eta_e(s) = (e,e)$. Now suppose $\eta_e(s) =\big( [s,s, e],[e,s,s]\big) = (e,e)$.  Then choosing some arbitrary $s' \in S$ we observe that 
$$[s,s,s'] = [s,s, [e,e,s']] = [[s,s, e], e ,s'] = [e,e,s'] = s'\,.$$
Thus, under our assumption, it must be the case that $s \in B_L$. Similarly,
$$[s', s,s] = [[s', e,e ] , s, s] = [s',e,[e, s,s ]] = [s', e,e ] = s'\,.$$
Thus, under our assumption, it must be the case that $s \in B_R$.  Together, we observe that $s \in B_L \cap B_R =: B$.  Thus we have  established that $B = \eta^{-1}_e(e,e)$ for any $e \in B$. \par 
Note that using Lemma \ref{lem:BClosed}, we see that $R_{e_1 e_2}(e_3) = [e_3, e_1, e_2] \in B$ and $L_{e_1, e_2}(e_3) = [e_1, e_2, e_3] \in B$ for all $e_1,e_2 , e_3 \in B$. That is, $B$ is an invariant set under the group of right and left translations generated by biunits (see Lemma \ref{lem:RLDiffs}).   Moreover, Lemma \ref{lem:BTrans} tells us that the group actions of right and left translations generated by biunits restricted to $B$ are transitive. Directly, we deduce
$$\big(R_{e_1 e} \times L_{e e_1} \big)\big(\eta_{e_1}(s) \big) = \big([[s,s,e_1], e_1,e], [e , e_1,[e_1 , s,s ]] \big) = \big([s,s,e], [e,s,s] \big) = \eta_e(s)\,.$$
\end{proof}
We are now in a position to state and prove the first theorem of this work.   
\begin{theorem}[Biunit Embedding]\label{thm:BiunitEmbed}
Let $S$ be a Lie semiheap. Then the subset of biunits $B := B_L \cap B_R$ is a closed embedded Lie subheap $\iota : B \hookrightarrow S$.  
\end{theorem}
\begin{proof}
If $B = \emptyset$, $B= S$, or $B = \{ e\}$, then the theorem is trivially satisfied. From this point on, we disregard these cases and consider biunital Lie semiheaps with multiple biunits.  \par  
As we have a biunital Lie semiheap, we can choose some $e \in B_L \cap B_R$. A standard result from semiheap theory tells us that the retract with respect to $e$, defined as $s_1 \cdot_e s_2 := [s_1, e, s_2]$, produces an associative binary product such that $e$ acts as the identity element. \par
We then define the maps
\begin{align*}
& \psi_L =  \Id_S : S \rightarrow S\,,  && \psi_R = \Id_S : S \rightarrow S \\
& s \mapsto s \cdot_e e = s \,,  && s \mapsto e \cdot_e s =s,
\end{align*}
which are translations evaluated at $e$, and so identity maps.   Note that the differential maps $(\rmd \psi_L)_{e}$ and  $(\rmd \psi_R)_{e}$ are the identity maps on $\sT_{e}S$. Then applying the Inverse Function Theorem (see \cite[Theorem 4.5]{Lee:2013}) to the map $(s_1, s_2) \mapsto (s_1, s_1 \cdot_e s_2)$, we know there exists an open neighbourhood $U \subset S$ of $e$ equipped with a smooth local inversion map, rendering $(U, \cdot_e)$ into a local Lie group.\par 
Next, consider $B_U := B \cap U$, which, due to Corollary \ref{cor:BClosedSubSet}, we know is a closed subset of the local Lie group.  Furthermore, Lemma \ref{lem:BClosed} ensures that, provided $e_1, e_2 \in B_U$ are close enough to $e$, then $e_1 \cdot_e e_2 = [e_1, e, e_2]\in B_U$. Similarly, the local inverse is given by $e^{-1}_1 :=  [e, e_1, e]$ (as standard). Lemma \ref{lem:BClosed} ensures that, provided $e_1$ is close enough to $e$, $e^{-1}_1 \in B_U$.  Thus, $(B_U, \cdot_e)$ is a closed local subgroup of the local Lie group $(U, \cdot_e)$. \par 
By Cartan's Closed Subgroup Theorem applied to local Lie groups (see  Lee \cite[Theorem 20.12]{Lee:2013}, whose proof naturally covers the local case), $B_U$ is an embedded smooth submanifold of $U$. As this property holds for all $e \in B$, it must be the case that $B$ is locally an embedded submanifold everywhere. Therefore, $B$ is a closed, embedded smooth submanifold of $S$.  Then via  Lemma \ref{lem:BClosed}, we conclude that $\iota : B \hookrightarrow S$ is a closed embedded Lie subheap.
\end{proof}
As $B$ is a Lie heap, we may pick any $e \in B$ and form the retract at $e$, which is a Lie group $G_e := (B, \cdot_e)$, provided $B \neq \emptyset$. Moreover, all such Lie groups are canonically isomorphic and so (up to isomorphism) we can speak of `the' \emph{biunital retract Lie group} of a biunital Lie semiheap. Thus, we also have a Lie algebra (up to isomorphism) associated with any biunital Lie semiheap. 
%
%
\subsection{The Infinitesimal Theory of Biunital Lie Semiheaps}\label{subsec:InfTheBiuLieSemiHeaps}
We will want to have a handle on infinitesimal transformations, and so the notion of a one-parameter group requires generalisation.  The approach we take requires our Lie semiheaps to be biunital. In particular, while mimicking one-parameter groups, we wish to avoid the definitions collapsing to heaps and thus groups.  First, we note that $I := [0, \epsilon) \subset \R_{\geq 0}$ comes with the smooth local semiheap structure $[t_1, t_2, t_3] := t_1 + t_2 + t_3$.  We remark that we do not obtain a local Lie semiheap by our definition, as $I$ is a manifold with boundary.
\begin{definition}\label{def:OneParamSemiHeap}
Let $S$ be a biunital Lie semiheap. A \emph{one-parameter semiheap} is a smooth forward-pointing curve
$$\varphi : I \rightarrow S\,,$$
with initial point  and $\varphi(0) = e \in B$, that is a local semiheap homomorphism, where $I := [0, \epsilon) \subset \R_{\geq 0}$ is equipped with the ``$+++$'' local  semiheap structure.
\end{definition}
To spell this out,  we have $I = [0, \epsilon)$, and given $t_1,t_2,t_3 \in I$, if $t_1+t_2+t_3 \in I$, then $\varphi(t_1+t_2+t_3) = [\varphi(t_1),\varphi(t_2),\varphi(t_3)]$.  Note that the image of a one-parameter semiheap is totally commutative, i.e., the $6$ permutations of the ordering of the elements are all equal. \par 
Note that $\varphi(t) =  \varphi(0+ t + 0) = [e, \varphi(t), e]$, with $\varphi(0) = e$. We define the map $\mathsf{I}_e : S \rightarrow S$ as $\mathsf{I}_e(s) := [e, s,e]$, noting that  $\mathsf{I}_e\big( \mathsf{I}_e(s)\big) = s$.  That is, $\mathsf{I}_e$ is a smooth involution on $S$. We then observe that the image of any one-parameter semiheap must lie entirely within the fixed-point set  of $\mathsf{I}_e$, that is
$$\mathrm{Im}(\varphi) \subseteq \mathrm{Fix}(\mathsf{I}_e) := \big \{ s \in S~~|~~ [e,s,e] = s \big\} \,.$$
As we have a smooth involution, $(\rmd \mathsf{I}_e)_e : \sT_e S \rightarrow \sT_e S$ is a linear involution.  Consequently, the tangent space at $e \in B$ splits into eigenspaces of $(\rmd \mathsf{I}_e)_e $;
\begin{equation}
\sT_e S  = V_{-1} \oplus V_{+1}\,.
\end{equation}
Differentiating $\varphi(t) = \mathsf{I}_e\big(\varphi(t) \big)$ at $t = 0^+$ and setting $\mathsf{v}:= \dot{\varphi}(0^+) \in \sT_e S$, we obtain 
\begin{equation}
\mathsf{v} = (\rmd \mathsf{I}_e)_e(\mathsf{v})\,.
\end{equation} 
Thus, the initial velocity of a non-trivial one-parameter semiheap must be an element of $V_{+1}$. \par 
Suppose  $e_1 \in B$; then $\mathsf{I}_e(e_1) = [e, e_1, e] = e^{-1}_1$, where the inverse is understood as being in the retract Lie group $G_e := (B, \cdot_e)$. As the derivatives of group inversions at the identity is $- \Id_{\sT_e B}$, we observe that $\sT_e B \cong V_{-1}$. This implies that any one-parameter semiheap whose image lies in $B$, or equivalently, has $\mathsf{v} \in \sT_eB$, must satisfy $\mathsf{v} \in V_{-1} \cap V_{+1}$, and thus $\mathsf{v} =0$, and so the one-parameter semiheap is trivial, i.e., $\varphi(t) = e \in B$.  \par 
We note that $\mathsf{I}_e$  acts as an anti-automorphism  of the local Lie group $(U, \cdot_e)$, i.e., 
$$\mathsf{I}_e(s_1 \cdot_e s_2) = [e, [s_1, e, s_2], e]= [e, s_2, [e, s_1, e] = [[e, s_2, e ], e ,[e, s_1,e]] = \mathsf{I}_e(s_2) \cdot \mathsf{I}_e(s_1)\,.$$
This observation is vital in constructing Lie triple system algebroids, and simply having an almost-product structure is insufficient. 
\begin{remark}
The opposite Lie semiheap $S^\textrm{op}$ is defined as $(S, [-,-,-]^\textrm{op})$ with $[s_1,s_2,s_3]^\textrm{op} := [s_3, s_2, s_1]$.  In general, a non-abelian Lie semiheap and its opposite are \ul{not} isomorphic. However, for biunital Lie semiheaps we have $\mathsf{I}_e \big([s_1,s_2, s_3] \big) = [\mathsf{I}_e(s_1), \mathsf{I}_e(s_2), \mathsf{I}_e(s_3)]^\mathrm{op}$.  Thus, if there is at least one biunit, then $S \stackrel{\sim}{\rightarrow} S^\mathrm{op}$.
\end{remark}
\begin{proposition}
Let $S$ be a proper biunital Lie semiheap. Then non-trivial one-parameter semiheaps $\varphi:  I \rightarrow S$ (see Definition \ref{def:OneParamSemiHeap}) always exist.
\end{proposition}
\begin{proof}
Let $S$  be a biunital Lie semiheap.  For any $e \in B$, we may select a non-zero $\mathsf{v} \in V_{+1}$ and consider the local one-parameter subgroup $\phi(t) = \exp_e(t \, \mathsf{v})$ of the local Lie group $(U, \cdot_e)$.  As $\mathsf{I}_e$ is an anti-automorphism fixing $e$ and the image of $\phi$ is abelian, $\mathsf{I}_e\big( \phi(t)\big)$ is a local one-parameter subgroup with initial velocity $(\rmd \mathsf{I}_e)_e (\mathsf{v})=\mathsf{v}$.  Uniqueness of one-parameter groups implies that $\mathsf{I}_e \big(\phi(t) \big) = \phi(t)$, and thus  $\mathrm{Im}(\phi) \subseteq \mathrm{Fix}(\mathrm{I}_e)$ as required by a one-parameter semiheap. \par 
Restricting attention to  $t_1, t_2, t_3 \in  I := [0, \epsilon)$, provided $t_1 + t_2 + t_3 \in I$, we calculate that
\begin{align*}
\phi(t_1 + t_2 + t_3) &= \phi(t_1) \cdot_e \phi(t_2) \cdot \phi(t_3) = [[\phi(t_1), e, \phi(t_2)], e, \phi(t_3)]\\
&=[\phi(t_1), [e, \phi(t_2), e], \phi(t_3)] =  [\phi(t_1), \mathrm{I}_e(\phi(t_2)), \phi(t_3)]\\
&= [\phi(t_1), \phi(t_2), \phi(t_3)]\,.
\end{align*} 
Thus, every $\mathsf{v} \in V_{+1}$ generates a one-parameter semiheap by  defining $\varphi := \phi|_I$.\par 
We know $\sT_e B \cong V_{-1}$.  This implies that $\dim(V_{+1})= \dim(\sT_eS) - \dim(\sT_e B)$. For proper biunital semiheaps (via Theorem \ref{thm:BiunitEmbed}) $\sT_e B$ is a proper subspace of $\sT_e S$.  This means that $\dim(V_{+1}) > 0$ and so we may always choose a non-zero $\mathsf{v} \in V_{+1}$. \par 
Thus, in conclusion,  non-trivial one-parameter semiheaps always exist if $S$ is a proper biunital Lie semiheap.
\end{proof}
As $B$ is a smoothly embedded submanifold of $S$, we may consider the restricted tangent bundle $\sT S|_B \rightarrow B$. For every $e \in B$, $\mathrm{I}_e(s) := [e,s,e] $ is a smooth involution  of $S$. The differential at the biunit $\big( \rmd \mathrm{I}_e\big)_e :\sT_e S \rightarrow \sT_e S$ is a linear involution and thus defines, via its eigenspaces,  the decomposition $\sT_e S \cong V_{-1}(e) \oplus V_{+1}(e) $.  Recall that $V_{-1}(e) \cong \sT_e B$. The dimension of $T_eB$ (and therefore $V_{-1}(e)$) is locally constant on connected components of $B$. Consequently, the complementary subspace dimension, $\dim(V_{+1}(e)) = \dim(T_eS) - \dim(V_{-1}(e))$, is  locally constant. Furthermore, right and left translations by biunits act transitively on $B$, which implies that $\dim(V_{+1}(e)) =k$ is globally constant across all of $B$. Because the base space $B$ is an embedded submanifold and each vector space $V_{+1}(e)$ is defined by a smoothly varying linear involution with constant eigenspace dimensions, the disjoint union $\coprod_{e \in B} V_{+1}(e)$ naturally forms a smooth vector subbundle of $TS|_B$. In short, associated with any proper biunital Lie semiheap is the vector bundle $\pi :A \rightarrow B$ with the fibres being given by the spaces $V_{+1}(e) \stackrel{\sim}{\rightarrow} \R^k$. Note that as $B$ is diffeomorphic to a Lie group, its retract Lie group, $\pi :A  \rightarrow B$ is trivialisable, but not canonically so. Defining a trivialisation is equivalent to selecting a particular $e_0 \in B$ to construct the retract Lie group. \par  
Recall that there always exists an open neighbourhood of any $e \in B$ such that $(U, \cdot_e)$ is a local Lie group, and that $\mathrm{I_e}$ acts as a local Lie group anti-homomorphism. Thus, at the level of the Lie algebra, $\theta := \big(\rmd \mathrm{I}_e\big)_e$ acts as a Lie algebra anti-homomorphism.   Thus, for $\mathsf{u}, \mathsf{v}\in V_{+1}(e)$, noting $\theta(\mathsf{u}) = \mathsf{u}$  and $\theta(\mathsf{v}) = \mathsf{v}$, we observe that $\theta([\mathsf{u}, \mathsf{v}]) = - [\mathsf{u}, \mathsf{v}]$, and so $[\mathsf{u}, \mathsf{v}] \in \sT_e B$. Thus, we have established that $\Sec(A)$ is \underline{not} closed under the standard Lie bracket of vector fields. In other words, $\pi : A \rightarrow B$ is not directly a Lie algebroid. For details of the general theory of Lie algebroids, the reader may consult  Mackenzie \cite{Mackenzie:2005}.\par    
However, $\Sec(A)$ is closed under the ternary structure $[\zx, \eta, \zeta] := [[\zx, \eta],\zeta]$, which defines a Lie triple system on the space of sections of $A$ (see Jacobson  \cite{Jacobson:1949} for Lie triple systems). To see the closure, we observe that 
$$\theta([\mathsf{u}, \mathsf{v}, \mathsf{w}] ) = \theta([[\mathsf{u}, \mathsf{v}] , \mathsf{w}] ) =- ([\theta [\mathsf{u}, \mathsf{v}] , \mathsf{w}]) =  [[\mathsf{u}, \mathsf{v}], \mathsf{w}] = [\mathsf{u}, \mathsf{v}, \mathsf{w}]\,,$$
for all $\mathsf{u}, \mathsf{v}, \mathsf{w}\in V_{+1}(e)$.
\begin{remark}
For every $e \in B$, $\sT_e S \cong V_{-1}(e) \oplus V_{+1}(e)$ is a  $\Z_2$-graded Lie algebra , i.e., 
$[V_{-1}(e) , V_{-1}(e)] \subseteq V_{-1}(e)$,  $[V_{-1}(e) , V_{+1}(e)] \subseteq V_{+1}(e)$ and $[V_{+1}(e) , V_{+1}(e)] \subseteq V_{-1}(e)$. Thus, the even subspace is $V_{-1}(e)$ and the odd subspace is $V_{+1}(e)$. This seems unintuitive, but recall $\theta$  is a Lie algebra anti-homomorphism. (We do not have a Lie superalgebra.)
\end{remark}
We are thus led to the second theorem of this work.
\begin{theorem}[Lie Triple Algebroid]\label{the:LieTripSys}
Let $S$ be a proper biunital Lie semiheap, and let $B$ denote the embedded submanifold of biunits. Then the vector bundle $\pi :A \rightarrow B$ defined above is a Lie triple algebroid (see Definition \ref{def:LieTripSymAlg}). 
\end{theorem}
\begin{proof}
In the preceding discussion we have established that $\Sec(A)$ is closed under the ternary structure  $[\zx, \eta, \zeta] := [[\zx, \eta],\zeta]$. It is a standard result that a Lie bracket defines a Lie triple system in exactly this way (see \cite{Jacobson:1949}).  Thus, we have 
\begin{itemize}
\itemsep0.5em
\item \emph{Skew-symmetric Property}: $[\zx,  \eta, \zeta] = - [\eta, \zx, \zeta]$;
\item \emph{Cyclic Identity}: $[\zx, \eta, \zeta] + [\zeta, \zx, \eta ] + [\eta, \zeta, \zx]=0$;
\item \emph{Fundamental Identity:} $[\zx, \eta, [\zeta, \kappa, \lambda]] =  [[\zx, \eta, \zeta], \kappa, \lambda] + [\zeta, [\zx, \eta, \kappa], \lambda] + [\zeta, \kappa, [\zx, \eta, \lambda]]$.
\end{itemize} 
We claim the anchor map $\rho : \Sec(A) \times \Sec(A) \rightarrow \Vect(B)$ is given by $\rho_{\zx \eta}:= [\zx, \eta]$, which is clearly antisymmetric. A direct calculation establishes this, i.e.,  for any $f \in C^\infty(B)$
 $$[\zx, \eta , f \, \zeta] = [[\zx, \eta], f \, \zeta] = [\zx, \eta](f) \, \zeta + f \, [[\zx, \eta], \zeta] = \rho_{\zx \eta}(f) \, \zeta + f \, [\zx, \eta , \zeta]\,. $$
The tensorial property  $[f \, \zx , \eta, \zeta] = [\zx, f\, \eta, \zeta] = f \, [\zx, \eta, \zeta]$ follows from direct calculation.  Note that $V_{+1}(e)$ is transverse to $V_{-1}(e) \cong \sT_e B$ for all $e \in B$. Thus, $\zx(f) =0$  for all $\zx \in \Sec(A)$ and  $f \in C^\infty(B)$.  This means that $[f\, \zx, \eta] = [\zx, f \, \eta] = f \, [\zx, \eta]$.  Hence, 
$$[f \zx , \eta, \zeta] = [\zx, f\,\eta, \zeta] =[f[\zx, \eta], \zeta] = f \, [[\zx, \eta], \zeta] =f\, [\zx, \eta, \zeta]\,. $$
\end{proof}
\smallskip
\noindent \textbf{Comments:}\
\begin{enumerate}[i)]
\item By construction we have $A \cong N(B) \subset \sT S|_B$, that is, $A$ can be identified with the normal bundle $N(B)$ to the submanifold of biunits. We then have a short exact sequence of vector bundles over $B$:
$$0 \rightarrow \sT B \rightarrow \sT S |_B \xrightarrow{\;\tau\;} A \rightarrow 0\,.$$
The inclusion  map $\sigma_e : A_e \hookrightarrow \sT_e S$, due to the decomposition  $\sT_e S = V_{-1}(e) \oplus V_{+1}(e)$ for all $e \in B$, defines a canonical splitting of the short exact sequence.  That is, we have a canonical right-splitting section $\sigma :  A \rightarrow \sT S|_B$.  Moreover, the anchor map $\rho : \Sec(A) \times \Sec(A) \rightarrow \Vect(B)$ measures the failure of the distribution $A \subset \sT S|_B$ to be integrable in the sense of Frobenius.
\item We have the isomorphism of $C^\infty(B)$-modules $\Sec(\sT S|_B) \stackrel{\sim}{\rightarrow} \big( C^\infty(B) \otimes_R\ \mathfrak{g}\big) \oplus \Sec(A)$, where $\mathfrak{g}$ is the Lie algebra associated with the retract Lie group of $B$ (which is defined up to isomorphism).
\item The sections $\Sec(A)$  represent the directions transverse to the submaifold of biunits. Thus, one may consider Lie triple algebroids as the analogue of Lie algebroids in the setting of Lie semiheaps.
\item If $B = \{ e\}$, then the anchor map must vanish, which implies $[\zx , \eta] =0$, which in turn implies that $[\zx, \eta , \zeta] =0$ for all $\zx, \eta$ and $\zeta \in \Sec(A)$.  That is, the Lie triple algebroid has the trivial (or zero) structure. 
\item  Although $A$ is always trivialisable,  in general, there is no canonical or preferred trivialisation; one has to select some $e_0 \in B$ to construct the retract Lie group of $B$, and this we view as the `affine nature' of the construction.
\end{enumerate}
\begin{example}\label{exa:nxnMatrices}
Consider the smooth manifold $S = M_n(\R) \cong \R^{n^2}$, which we equip with the ternary operation $[X, Y, Z] := X\, Y^T Z$. As matrix multiplication and transposition are smooth operations, we obtain a  Lie semiheap. The submanifold of biunits is given by the set of orthogonal matrices, i.e.,  $B = O_n(\R)$.  Thus, $S$ is a proper biunital Lie semiheap.  \par 
At a biunit $E \in O_n(\R)$, $\mathsf{I}_ E(X) = E \, X^{T} E $. As $\sT_E M_n(\R) \cong M_n(\R)$, we have 
$$\big( \rmd \mathsf{I}_E \big)_E (V) = E \, V^T E\,.$$
The tangent bundle decomposes as $\sT_E S = V_{-1}(E) \oplus V_{+1}(E)$, where
$$V_{-1}(E) = \big \{ V \in M_n(\R)~~|~~ E \, V^T E = - V \big\}\,, \qquad \text{and} \quad V_{+1}(E) = \big \{ V \in M_n(\R)~~|~~ E \, V^T E =  V \big\}\,.$$
For $V \in V_{-1}(E)$, $V^T E = - E^T V$, which implies $(E^T V)^T = - E^T V$. Thus, $E^TV$ is a skew-symmetric matrix.  This directly leads to $V_{-1}(E) = E \cdot \mathfrak{o}(n) \cong \mathfrak{o}(n) = \sT_E B$ as expected. Similarly,  for $V \in V_{+1}(E)$, $V^T E = E^T V$, which implies $(E^T V)^T =  E^T V$. Thus, $E^TV$ is a symmetric matrix.  Thus, we obtain the well-known linear decomposition of $M_n(\mathbb{R})$ into symmetric and skew-symmetric subspaces. This directly leads to $V_{+1}(E) = E \cdot \mathrm{Sym}_{n}(\R) \cong \mathrm{Sym}_{n}(\R)$. \par 
We observe that the vector bundle $\pi : A \rightarrow O_n(\R)$ is trivialisable and  
$$A \cong O_n(\R) \times \mathrm{Sym}_{n}(\R)\,.$$ 
Thus, $\Sec(A) \cong C^\infty\big( O_n(\R), \mathrm{Sym}_n(\R)\big)$. That is, sections of $A$ are identified with  $n\times n$ symmetric matrices with entries being smooth functions on the Lie group $O_n(\R)$.  The Lie triple algebroid is given by
$$[\zx, \eta, \zeta] := [[\zx,\eta],\zeta]\,, \quad \text{and} \quad \rho_{\zx \eta} := [\zx, \eta]\,,$$
where the binary bracket is the standard matrix commutator of matrices with functions as entries. Note the abuse of notation here, as written, sections of $A$ are matrix-objects rather than derivations of $C^\infty(O_n(\R))$. To recover a tangent vector at a point $E \in O_n(\R)$, we undo the earlier isomorphism and set $X_\zx(E) := E \cdot \zx(E) \in \sT_E M_n(\R)$.    
\end{example}
\begin{example}[The conformal causal Lie semiheap] \label{exa:ConfCausal}As a physically relevant specific case of Example \ref{ema:OlSemi}, consider $G = \mathrm{SU}(2,2)$, which is the double cover of the identity connected component of the conformal Lorentz group. This group consists of $4\times 4$ complex matrices $g$, such that $\det(g) = 1$. Defining $\mathrm{J} = \textrm{diag}(\mathrm{I}_2 , - \mathrm{I}_2)$,  the associated involution is given by $\Theta(g) =  J (g^\dagger)^{-1} J$.  Here $K \cong \mathrm{S}(\mathrm{U}(2) \times \mathrm{U}(2))$. Then, at the level of the Lie algebra, the involution induces a decomposition $\mathfrak{su}(2,2) = \mathfrak{k}\oplus \mathfrak{p}$. The compact component generates pure spatial rotations and a compact ``global time'' translation. The non-compact component generates standard relativistic boosts, scale transformations, spatial translations, and accelerating transformations.  Elements of the causal subspace $\mathfrak{p}_c \subset \mathfrak{p}$ can be parametrised using spacetime vectors $(x^0, x^1, x^2, x^3) \in \mathbb{R}^{1,3}$. The elements of  $\mathfrak{p}_c$ are $4 \times 4$ block off-diagonal matrices of the form
$$ X = \left(\begin{matrix}0&v\\ v^\dagger &0\end{matrix}\right)\,.$$
The space $\mathfrak{p}_c$ is isomorphic to the space of $2 \times 2$ Hermitian matrices via its off-diagonal blocks. Thus any element $X \in \mathfrak{p}_c$  can be represented as a matrix 
$$v = \begin{pmatrix} x^0 + x^3 & x^1 - \rmi \,x^2 \\ x^1 + \rmi\, x^2 & x^0 - x^3 \end{pmatrix}\,.$$
The  Minkowski metric $\rmd s^2 = (x^0)^2 - (x^1)^2 - (x^2)^2 - (x^3)^2$ is directly encoded by the determinant, i.e., $\det(v) = (x^0)^2 - (x^1)^2 - (x^2)^2 - (x^3)^2$. The open, convex, pointed cone $C \subset \mathfrak{p}_c$ is defined by requiring the block $v$ to be strictly positive definite. In terms of vector components, this translates to the causal conditions:
$$x^0 > 0 \quad \text{and} \quad \det(v) > 0\,.$$
That is, we have an open future-directed light cone. This cone is invariant under the adjoint action of the maximal compact subgroup $K$, i.e.,  $\mathrm{Ad}(K)C = C$.   The open Ol'shanskii semigroup $S^\textnormal{o}$ is defined as the image of $\mathrm{SU}(2,2) \times C \rightarrow \mathrm{SU}(2,2)_\C$ of the map $(g, X) \mapsto g \, \exp(\rmi\, X)$.  We remark that as a real manifold, $S^\textnormal{o}$ is $19$-dimensional, and $\mathrm{SU}(2,2)_\C$ is $30$-dimensional.  We then define the Lie semiheap operation $[s_1, s_2, s_3] := s_1 s_2^* s_3$ (see Example \ref{ema:OlSemi}), noting that the submanifold of biunits   is $B \stackrel{\sim}{\rightarrow} \mathrm{S}(\mathrm{U}(2) \times \mathrm{U}(2)) \times C$, using the polar decomposition. Thus, we obtain a proper biunital Lie semiheap. \par 
To describe the associated Lie triple algebroid, we proceed as follows.  Using the parametrisation of $S^\textnormal{o}$ and the Cartan decomposition, we observe that for any $s \in S^\textnormal{o}$, $\sT_e  S^\textnormal{o} \stackrel{\sim}{\rightarrow} \mathfrak{su}(2,2) \oplus \mathfrak{p}_c \cong (\mathfrak{k} \oplus \mathfrak{p}_c) \oplus V_{+1}(e)$. Thus, $V_{+1}(e) \cong \mathfrak{su}(2,2)\setminus \mathfrak{k} = \mathfrak{p}$. We then have the non-canonical isomorphism $A \stackrel{\sim}{\rightarrow} \big( \mathrm{S}(\mathrm{U}(2) \times \mathrm{U}(2)) \times C\big)\times \mathfrak{p}$. Clearly, $\mathfrak{p}$ is not closed under the standard Lie bracket, but it is under the nested Lie bracket. Thus, we obtain a Lie triple algebroid.  \par
The physical interpretation of a biunit is a complete observer state. In particular, using the polar decomposition, we have $e = (k, v) \in B$. The element $k \in \mathrm{S}(\mathrm{U}(2) \times \mathrm{U}(2))$ represents a spatial orientation and a global time phase. The element $v \in C$ is a choice of a future-directed timelike vector.  Thus, a biunit defines a complete observer reference frame; in other words, $B$ is identified with the state space of observer reference frames. In the standard language of involutive semigroups, observers correspond to unitary elements. Sections of the associated Lie triple algebroid are interpreted as complete observer-dependent non-compact infinitesimal conformal transformations.
\end{example}
\begin{remark}
Using $\mathrm{SU}(1,1)$ in place of $\mathrm{SU}(2,2)$ and $1+1$-dimensional Minkowski spacetime, in Example \ref{exa:ConfCausal}, we obtain Howe's oscillator (open) semigroup (see \cite{Howe:1988}) and the related proper biunital Lie semiheap.  In this case $B \cong \mathrm{U}(1) \times C$, where $k \in \mathrm{U}(1)$ is the observers global time phase.\end{remark}
Let $S$ be a proper biunital semiheap, and $\pi : A \rightarrow B$ be its associated Lie triple algebroid.  We treat the vector bundle abstractly and will not fix a trivialisation. Let us select a (non-zero) section $\zx \in \Sec(A)$ and a biunit $e \in B$. Then $\mathsf{v} = \zx(e) \in V_{+1}(e)$ generates a local one-parameter semiheap (see Definition \ref{def:OneParamSemiHeap}) $\varphi_e : I \rightarrow S$ that satisfies $\varphi_e(0)= e$ and $\dot{\varphi}_e(0^+) = \mathsf{v}$.  We then define a \emph{left-derived} vector field $X_\zx \in \Sec(S)$ via
\begin{equation}
X_\zx(s) := \left . \frac{\rmd}{\rmd t}\right|_{t = 0^+} \, L_{e \varphi_e(t)}(s)\,,
\end{equation}
where $L_{e \varphi_e(t)}(-) :=  [e, \varphi_e(t), -]$ is a smooth family of left translations (see Definition \ref{def:RLActions}). The notion of a \emph{right-derived} vector field follows similarly using right translations. 
\begin{proposition}\label{prop:XCovariant}
Let $S$ be a proper biunital Lie semiheap and $\pi :A \rightarrow B$ be its associated Lie triple algebroid. Furthermore, let $\zx \in \Sec(A)$ be an arbitrary (non-zero) section.  Then the associated left-derived vector field $X_\zx \in \Vect(S)$  transforms covariantly under a change of biunit $e \mapsto e'$. 
\end{proposition}
\begin{proof}
Let $e \in B$ be the chosen biunit that is used to define $X_\zx$ in the above construction.  Suppose we select some other biunit $e' \in B$.  We note that $e' = [e',e,e] = L_{e'e}e$. We now consider the one-parameter semiheap defined as $\varphi_{e'}(t) = [e', e,\varphi_e(t)] = L_{e'e}\varphi_e(t)$, which satisfies $\varphi_e'(0)= e'$. Setting $\Phi :=  L_{e'e}$ for convenience, we note that, via direct calculation  (which we omit), $L_{e'\varphi_{e'}(t)} s = \Phi \big(L_{e \varphi_e(t)}(\Phi^{-1}s) \big)$. Taking the derivative, we arrive at 
$$X'_\zx(s) := \left . \frac{\rmd }{\rmd t}\right|_{t = 0^+} L_{e' \varphi_{e'}(t)} = \big ( \rmd \Phi\big)_{\Phi^{-1}(s)} \left( X_\zx \big(\Phi^{-1}(s)\big) \right)\,.$$
Rephrasing this slightly, we have $X'_\zx = \Phi_* X_\zx$. 
\end{proof}
The above proposition can be restated as the following diagram is commutative.
\renewcommand\arraystretch{1.5}
$$\begin{array}{ccc} \Sec(A) & \xrightarrow{\quad X\quad} & \Vect(S) \\ \parallel & & \big\downarrow \Phi_* \\ \Sec(A) & \xrightarrow{\quad X'\quad} & \Vect(S) \end{array}$$
It is important to note that the constructed maps are $\R$-linear rather than module homomorphisms.  
\begin{proposition}
Let $S$ be a proper biunital Lie semiheap and $\pi :A \rightarrow B$ be its associated Lie triple algebroid.  Then the map  $X: \Sec(A) \rightarrow \Vect(S)$  satisfies
$$ [[X_\zx, X_\eta], X_\zeta] = X_{[[\zx, \eta], \zeta]} = X_{[\zx, \eta, \zeta]}\,,$$
for all $\zx, \eta, \zeta \in \Sec(A)$. 
\end{proposition}
\begin{proof}
Let $e \in B$ be the biunit used to construct the map $X: \Sec(A) \rightarrow \Vect(S)$. Furthermore, let $G_e = (U, \cdot_e)$ be the local Lie group around the bunit. We then define a left (local) Lie group action  $\sigma : G_e \times S \rightarrow S$  by $\sigma(g, s) := g \cdot_e s = [g , e, s]$ (one can directly check we have an action).  We then observe that $[e, \varphi_\zx(t), s] = \varphi_\zx(t) \cdot_e s = \sigma(\varphi_\zx(t), s)$ using para-associativity.  Thus, $X_\zx$ is the fundamental vector field on $S$ corresponding to $\zx(e) \in \mathfrak{g}_e$ under the left action.  Thus, as standard we have $[X_\zx, X_\eta] = - X_{[\zx, \eta]}$.   However, we know $[\zx,\eta ]$ is not a section of $A$. Taking the bracket again we obtain $ [[X_\zx, X_\eta], X_\zeta] = X_{[[\zx, \eta], \zeta]} = X_{[\zx, \eta, \zeta]}$.
\end{proof}
By defining $[X_\zx, X_\eta, X_\zeta] :=  [[X_\zx, X_\eta], X_\zeta]$, we view $X$ as a  Lie triple system homomorphism. Moreover, using Proposition  \ref{prop:XCovariant}, together with the fact that pushforwards with respect to diffeomorphisms commute with the Lie bracket of vector fields, the Lie triple system realised in $\text{Vect}(S)$ is unique up to canonical Lie triple system isomorphism.  The vector fields $X_\zx \in \Vect(S)$ constructed from sections $\zx \in \Sec(A)$ are understood as ``fundamental vector fields''  that infinitesimally capture the Lie semiheap structure in `directions off' the submanifold of biunits $\iota :B\hookrightarrow S$. We can then examine Lie derivatives of tensors (and tensor densities) with respect to these vector fields. Importantly, as we have a Lie  triple system, we have 
$$\big[[\mathcal{L}_{X_\zx}, \mathcal{L}_{X_\eta} ],\mathcal{L}_{X_ \zeta}\big]= \mathcal{L}_{X_{[\zx, \eta, \zeta]}}\,,$$
rather than the standard binary closure.  However, care must be taken as the map $X : \Sec(A) \rightarrow \Vect(S)$ depends on the chosen biunit $e \in B$. \par 
Given a section $\zx \in \Sec(A)$, we define an associated distribution  $\mathcal{D}_\zx \subset \sT S$  (generally, in the sense of Stefan and Sussmann, this will be singular) as
\begin{equation}
\mathcal{D}_\zx(s) := \Span_{\R} \big \{X^e_\zx(s)~~|~~ e \in B  \}\,,
\end{equation} 
where $X^e : \Sec(A) \rightarrow \Vect(S)$ denotes the map with the required biunit made explicit. We then define a \emph{$\zx$-invariant function} as a function $f \in C^\infty(S)$, such that $\mathcal{L}_{X^e_\zx} f =0$ for all $e \in B$.  In other words, $f$ is a first-integral of $\mathcal{D}_\zx$, i.e., $\rmd f( \mathcal{D}_\zx)=0$. Similarly, we define a \emph{$\zx$-invariant vector field} as a vector field $Y \in \Vect(S)$ that is a strong symmetry of $\mathcal{D}_\zx$, i.e., $[Y, \mathcal{D}_\zx]=0$. \par 
Looking towards infinitesimal actions of Lie semiheaps on manifolds, we propose the following definition using the associated Lie triple algebroid. 
\begin{definition}
Let $S$ be a proper biunital Lie semiheap and $\pi : A \rightarrow B$ be its associated Lie triple algebroid. Furthermore, let $M$ be a smooth manifold. Then an \emph{action of a Lie triple algebroid on a manifold} is a continuous $\R$-linear map $\Phi :  \Sec(A) \rightarrow \Vect(M)$ that satisfies
$$\Phi([\zx, \eta, \zeta]) = [[\Phi(\zx), \Phi(\eta)],\Phi(\zeta)]\,,$$
for all sections $\zx, \eta$ and $\zeta \in \Sec(A)$.
\end{definition}
In other words, an action of a Lie triple algebroid on a manifold is a (continuous) homomorphism of Lie triple systems $\Phi: \Sec(A) \rightarrow  \Vect(M)$, where the triple bracket of vector fields on $M$ is the nested Lie bracket.
\begin{example}
Continuing Example \ref{exa:ConfCausal}, the associated proper biunital Lie semiheap is described by $B \stackrel{\sim}{\rightarrow} \mathrm{S}(\mathrm{U}(2) \times \mathrm{U}(2)) \times C$, and $A_{(k,v)} \stackrel{\sim}{\rightarrow} \mathfrak{p}$. Consider the action
$$\Phi : \Sec(A) \rightarrow \Vect(B \times \mathbb{M})\,,$$
where $\mathbb{M}$ is four-dimensional Minkowski spacetime.  Thus, $\Phi(\zx)$  for any $\zx \in \Sec(A)$ is an observer-dependent infinitesimal non-compact conformal transformation acting on the extended kinematic space (also known as the observer bundle) $B \times \mathbb{M}$.   As $\sT(B \times \mathbb{M}) \cong \sT B \times \sT \mathbb{M}$, we have the decomposition
$$\Phi(\zx)|_{(k, v, x)} = \Phi_B(\zx)|_{(k, v, x)} + \Phi_{\mathbb{M}}(\zx)|_{(k, v, x)}\,.$$
Thus, the transformations, in general, do not preserve the (trivial) bundle structure of the extended kinematic space. It is important to note that these transformations are  invertible.  In particular, given a section $\zx \in \Sec(A)$ the section $- \zx$ exists, and due to $\R$-linearity,  $\Phi(-\zx) = - \Phi(\zx)\in \Vect(B \times \mathbb{M})$. Thus, at the infinitesimal level we have local invertibility. However, integrating the flows to construct global transformations  forces the integral curves to flow in one-direction as dictated by the causal cone.
\end{example}
%
%
\section{Concluding Remarks}\label{subsec:ConRem}
In this work, we have investigated the differential geometry of proper biunital Lie semiheaps.  We remark that examples of proper biunital Lie semiheaps are common, including the space of $n\times n$ real matrices and (open) Ol'shanskii semigroups. The core results are that the set of biunits is a submanifold, there is an associated Lie triple systems algebroid over the biunits, and that we may extend the sections of this algebroid to vector fields on the Lie semiheap to form a Lie triple system (defined uniquely up to a canonical isomorphism). Finally, we showed that all the possible vector fields associated with a section can be used to construct a distribution and that this distribution allows for a notion of functions and vector fields invariant under the action/flow of a section. The idea generalises to differential forms and tensors, though we have not presented details here. \par 
\medskip 

%
\noindent\textbf{Potential Future Directions:}\
\begin{enumerate}[i)]
\item Further examples of proper biunital Lie semiheaps are highly desirable to aid the understanding of the general theory. Submanifolds of $M_n(\R)$ may provide non-trivial prototypes.  For instance, we may restrict attention to the open submanifold  $\mathrm{GL}_n(\R)\subset M_n(\R)$ of invertible matrices. Note that the semiheap operation  $[X, Y, Z] := X \, Y^T Z$ closes on invertible matrices and $B = O_n(\R)$. More generally, transpose-closed Lie subgroups (and their connected components) $G \subset \mathrm{GL}_n(\R)$ form proper biunital Lie semiheaps with $B = G \cap O_n(\R)$. More interesting from a physics perspective would be to develop further examples along the lines of the conformal causal Lie semiheap (Example \ref{exa:ConfCausal}) but now using non-semisimple Lie groups and an invariant cone in a suitable Lie algebra. 
\item Applications in geometric mechanics and control theory seem tractable. In particular, Lie semiheaps offer a framework for smooth non-invertible symmetries/transformations; for example, these may be useful in the modelling of dissipative systems.  Proper biunital  Lie semiheaps offer a way to understand infinitesimal non-invertible symmetries.   It is expected that Lie triple systems will play a vital role here, akin to standard Lie algebras and invertible symmetries.  For example, we expect that geometric mechanics following Tulczyjew \cite{Tulczyjew:1976,Tulczyjew:1977} can be developed using proper biunital Lie semiheaps and their associated Lie triple algebroids. This awaits proper exploration.
\item  It is well-known that Lie triple systems give rise to solutions of the  Classical Yang--Baxter Equation (CYBE) via relative Rota--Baxter operators. This strongly suggests that smooth, parameter-dependent solutions to the CYBE  may be constructed using proper biunital Lie semiheaps. In particular, a relative Rota--Baxter operator (of weight $0$) we understand as a $C^\infty(B)$-linear map $R : \Sec(A) \rightarrow \Sec(A)$ that satisfies 
$$[R(\zx), R(\eta), R(\zeta)] = R \big( [R(\zx), R(\eta), \zeta]+ [R(\zx), \eta, R(\zeta)] + [\zx, R(\eta), R(\zeta)]\big)\,.$$ 
The descendant Lie triple algebroid is defined as 
$$[\zx, \eta, \zeta]_R := [R(\zx), R(\eta), \zeta]+ [R(\zx), \eta, R(\zeta)] + [\zx, R(\eta), R(\zeta)]\,,  \quad \text{and} \quad  \rho^R_{\zx \eta} :=  \rho_{R(\zx)R(\eta)}\,.$$
Following the standard arguments for Lie triple systems and relative Rota--Baxter operators, it can be shown that we do indeed obtain a Lie triple algebroid in this way. Furthermore, note $R([\zx, \eta, \zeta]_R) = [R(\zx), R(\eta), R(\zeta)]$.  Details of constructing (smoothly parametrised)  solutions to the CYBE and explicit examples remain completely unexplored.
\end{enumerate}

%
%
%
\medskip 

\noindent \textbf{LLM Statement:} Google Gemini (version 3.6 Flash) was used to assist in the conceptualisation of proofs and examples. All results were verified by the author.   
%
%

%
%
\newpage
\appendix
\section{Lie Triple Algebroids}\label{App:LieTripSymAlg}
Lie triple algebroids we view as the oidification or horizontal categorification of Lie triple systems. Heuristically, the vector space structure underlying a Lie triple system is replaced with a vector bundle. Our definition is guided by Theorem \ref{the:LieTripSys} rather than category theory.  The notion of an abstract Lie triple algebroid is the following.
\begin{definition}\label{def:LieTripSymAlg}
A \emph{Lie triple algebroid} is a vector bundle $\pi: A \rightarrow B$ such that the module of sections $\Sec(A)$ is equipped with a ternary operation $[-,-,-] : \Sec(A)^{\times^{3}} \rightarrow \Sec(A)$ and an anchor map $\rho_{- -}: \Sec(A) \times \Sec(A) \rightarrow \Vect(B)$ that satisfies the following:\\
\begin{enumerate}[a)]
\itemsep1em
\item \emph{Skew-symmetric Property}: $[\zx,  \eta, \zeta] = - [\eta, \zx, \zeta]$;
\item \emph{Cyclic Identity}: $[\zx, \eta, \zeta] + [\zeta, \zx, \eta ] + [\eta, \zeta, \zx]=0$;
\item \emph{Fundamental Identity:} $[\zx, \eta, [\zeta, \kappa, \lambda]] =  [[\zx, \eta, \zeta], \kappa, \lambda] + [\zeta, [\zx, \eta, \kappa], \lambda] + [\zeta, \kappa, [\zx, \eta, \lambda]]$;
\item  \emph{Leibniz Rule:} $[\zx, \eta , f \, \zeta] = \rho_{\zx \eta}(f) \, \zeta + f \, [\zx, \eta , \zeta]$;
\item \emph{Tensorial Property}:  $[f \, \zx , \eta, \zeta] = [\zx, f\, \eta, \zeta] = f \, [\zx, \eta, \zeta]$,
\end{enumerate}
\smallskip
for all $\zx, \eta, \zeta, \kappa, \lambda \in \Sec(A)$ and $f \in C^\infty(B)$.
\end{definition}
\noindent \textbf{Comment.}  The related notions of Lie triple system bundles and Lie triple system algebroids have independently been proposed by Banerjee \&  Mukherjee \cite{Banerjee:2026}. While the Lie triple algebroids introduced here are Lie triple system bundles, the anchor map and the identities it satisfies are different. Their motivation is the tangent bundle of locally symmetric spaces, while the motivation of the current paper is the geometry of Lie semiheaps. Note that the anchor map of a Lie triple algebroid takes two sections as the input, while the anchor map of a Lie triple system algebroid requires only a single section. The reader should consult \cite[Definition 7.1]{Banerjee:2026} for details.\par
\smallskip
Using these axioms, it can be shown that
\begin{subequations}
\begin{align}
&\rho_{\zx \eta} = - \rho_{\eta  \zx}\,,\\
& [\rho_{\zx\eta}, \rho_{\zeta\kappa}] = \rho_{[\zx,\eta,\zeta]\kappa} + \rho_{\zeta[\zx,\eta,\kappa]}\,.
\end{align}
\end{subequations}
Unlike the specific examples of such structures associated with proper biunital Lie semiheaps, the vector bundle structure may be non-trivialisable.
\begin{warning}
Lie triple algebroids should not be confused with Filippov algebroids (of order $3$) despite the clear similarities; see Grabowski \& Marmo \cite{Grabowski:2000} for details of Filippov algebroids. 
\end{warning}
Note that the derived anchor relation mirrors the construction in a Lie triple system, where the left-multiplication operator $L(x,y)z := [x,y,z]$ acts as an inner derivation. The commutator of two such derivations satisfies $[L(x,y), L(u,v)] = L([x,y,u], v) + L(u, [x,y,v])$. Furthermore, if the vector bundle is over a single point, and so identified with a vector space, then 
the anchor map vanishes, and we recover the standard notion of a Lie triple system (see Jacobson  \cite{Jacobson:1949}).
\begin{example}
Let $\pi : A \rightarrow B$ be a totally intransitive Lie algebroid, i.e., the anchor map is the zero map. In other words, $A$  is a bundle of Lie algebras. Then $[\zx, \eta, \zeta] := [[\zx, \eta], \zeta]$ defines a Lie triple algebroid structure on $A$. In this case, we have a bundle of Lie triple systems. 
\end{example}
\begin{warning}
We do not obtain a Lie triple algebroid by nesting the Lie bracket on a Lie algebroid other than in the totally transitive case.
\end{warning}
%
%

\end{document}